\documentclass[11pt]{article}
\usepackage{amssymb}
\usepackage{amsmath}
\usepackage{amsthm}
\usepackage{latexsym}
\usepackage{amsfonts}
\usepackage{graphicx}
\usepackage{graphics}

\newcommand{\dis}{\displaystyle}
\theoremstyle{plain}
\newtheorem{thm}{Theorem}[section]   % Αρίθμηση συνεχόμενη (όχι κατά Θεώρημα, Λήμμα κ.λπ.)
\newtheorem{prop}[thm]{Proposition}
\newtheorem{lem}[thm]{Lemma}

\newtheorem{Def}[thm]{Definition}%[section]
\newtheorem{Defs}[thm]{Definitions}%[section]

\theoremstyle{definition}
\newtheorem{rem}[thm]{Remark}
\newtheorem{rems}[thm]{Remarks}
\newtheorem{exm}[thm]{Example}
\newtheorem{note}[thm]{Note}
\newtheorem*{Proof}{Proof}

\newcommand{\bbb}[1]{\mbox{\boldmath$#1$}}

\newcommand{\fa}{\forall}

\newcommand{\el}{\ell}

\newcommand{\ra}{\;\rightarrow\;}

\newcommand{\Ga} {{\varGamma}}
\newcommand{\Sig} {{\varSigma}}
\newcommand{\de}{\delta }

\newcommand{\e}{\varepsilon }

\newcommand{\la}{\lambda }
\newcommand{\mi}{\mu }

\newcommand{\si}{\sigma }

\newcommand{\R}{\mathbb{R}}

\newcommand{\N}{\mathbb{N}}

\newcommand{\ssum}{\sum\limits}

\newcommand{\ld}{\ldots}

\newcommand{\sm}{\smallsetminus}

\begin{document}
\title{\bf The Classical Theorems of Measure Theory in connection with the Statistical convergence and some remarks on Steinhaus' Theorem}
\author{Christos Papachristodoulos}
\date{}
\maketitle
\begin{abstract}
We study Steinhaus' theorem regarding statistical limits of measurable real valued functions and we examine the validity of the classical theorems of Measure Theory for statistical convergences.
\end{abstract} \medskip
 \medskip \\
Mathematics Subject Classification (2000) 28A20
\section{Introduction}\label{sec1}
\noindent

Statistical convergences was introduced by Zygmund in a monograph in 1935 (see [9], vol. II, p. 181) and continued by Steinhaus and Fast a few years later (\cite{8}, \cite{4}). Since then, several related papers have been published, mainly on applications and generalizations of this notion of convergence (\cite{1}, \cite{2}, \cite{5}, \cite{7}).  Our aim is to deal with some remarks regarding Steinhaus' theorem (Th. 1.6 below) and statistical convergences of sequences of measurable functions. In order to be more concrete, it is convenient to start with the framework that this paper is based on, as well as, the definitions and known results that we need.
\begin{enumerate}
\item[$\bullet$] Throughout the paper, for the sake of simplicity, we consider the space $([0,1),\Sig,\la)$, where $\Sig$ is the $\si$-algebra of Lebesgue measurable subsets of $[0,1)$ and $\la$ is the Lebesgue measure.
\end{enumerate}

We note that all the results, that we will see, hold for any finite measure space and some of them, that we point out by remarks, holds for arbitrary, measure space.
\begin{enumerate}
\item[$\bullet$] All functions are real valued measurable functions defined on $[0,1)$.
\item[$\bullet$] $\N$ is the set of positive integers.
\item[$\bullet$] For $A\subseteq\N$ by $d(A)$ we denote the density of $A$, that is,
\[
d(A)=\lim_{n\ra\infty}\frac{|\{k\in A:k\le n\}|}{n},
\]
if of course the above limit exists.
\item[$\bullet$] If $I\subseteq[0,1]$, by $\chi_I$ we denote the characteristic function of $I$,
\[
\chi_I(a)=\left\{\begin{array}{cl}
                   1, & \text{if} \ \ x\in I \\
                   0, & \text{otherwise}
                 \end{array}\right.
\]
\item[$\bullet$] $C^+_0=\{(\e_n)|\e_n>0$ for $n=1,2,\ld$ and $(\e_n)\ra0\}$.
\item[$\bullet$] The notations $(f_n)\overset{a.e.}{\longrightarrow}f$, $(f_n)\overset{\la}{\longrightarrow}f$, $(f_n)\xrightarrow[]{al-u}f$ means respectively that, the sequence of real values measurable functions $(f_n)$ converges almost everywhere, in measure, almost uniformly to $f$.
\item[$\bullet$] $L^0=L^0([0,1))$ is the space of real valued measurable functions defined on $[0,1)$, where as usual we consider $f=g$, if $f(x)=g(x),\la-a.e$.

\end{enumerate}
\begin{Defs}\label{Defs1.1}
Let $(a_n)$ be a sequence in $\R$ and $a\in\R$.
\begin{enumerate}
\item[(a)] We say that $(a_n)$ converges statistically to $a$ and we write $(a_n)\overset{st}{\longrightarrow}a$, if, $\fa\;\e>0$
    $d(\{n\in\N:|a_n-a|\ge\e\})=0$.
\item[(b)] We say that $(a_n)$ is statistically Cauchy and we write $(a_n)$ is $st-C$, if, \vspace*{-0.4cm}
\end{enumerate}
\[
\fa\;\e>0 \ \ \exists\;n_0\in\N:d(\{n\in\N:|a_n-a_{n_0}|\ge\e\})=0.
\]
\end{Defs}

Regarding statistical convergences of numerical sequences we will need the following known results.
\begin{prop}\label{prop1.2}
Let $(a_n)$ be a sequence in $\R$ and $a\in\R$.

(i)\; $(a_n)\overset{st}{\longrightarrow}a\;\Leftrightarrow\;\exists\;K=\{k_1<k_2<\ldots\}
\subseteq\N,d(K)=1:(a_{k_n})\ra a$.

(ii)\; $(a_n)$ converges statistically $\Leftrightarrow\;(a_n)$ is $st-C$.

(iii)\; If moreover $(a_n)$ is bounded, then
\[
(a_n)\overset{st}{\longrightarrow}a\;\Leftrightarrow\;\lim_{n\ra\infty}\frac{1}{n}
\sum^n_{k=1}a_k=a.
\]

(iv)\; If $a_n\ge a$ for $n=1,2,\ld$, then
\[
\lim_{n\ra\infty}\frac{1}{n}\sum^n_{k=1}a_k=a\;\Rightarrow\;(a_n)\overset{st}{\longrightarrow}a.
\]
\end{prop}

(For the proof see \cite{2}, \cite{5}, \cite{7}).
\begin{rem}\label{rem1.3}
Obviously, Definitions \ref{Defs1.1} extend to arbitrary metric space and in this case Proposition \ref{prop1.2}(i) is true, while Proposition \ref{prop1.2}(ii) holds whenever the metric space is complete.
\end{rem}
\begin{Def}\label{Def1.4} {\em (\cite{4})} A sequence $(a_n)$ in $\R$ is said to be measurable if, there exist an at most countable subset $A$ of $\R$ such that the density of the sets $\{n\in\N:a_n<a\}$ exists for all $a\in\R\sm A$.
\end{Def}
\begin{Def}\label{Def1.5} {\em (\cite{4})} Let $f_n$, $f:[0,1)\ra\R$, $n=1,2,\ld$ be measurable functions.

a) We say that the sequence $(f_n)$ converges statistically almost everywhere to $f$ and we write $(f_n)\xrightarrow[]{st-a.e.}f$, if
\[
f_n(x)\overset{st}{\longrightarrow}f(x),
\]
$\la$-almost everywhere for $x\in[0,1)$.

(b) We say that the sequence $(f_n)$ converges statistically in measure or asymptotic statistically to $f$ and we write $(f_n)\xrightarrow[]{st-\la}f$, if,
\[
\fa\;\e>0:\la([|f_n-f|\ge\e])\overset{st}{\longrightarrow}0.
\]
\end{Def}

With the above notations Stainhaus' theorem stands as follows.
\begin{thm}\label{thm1.6}
{\em (\cite{4})} (i) If $(f_n)\xrightarrow[]{st-a.e.}f\;\Rightarrow\;(f_n)\xrightarrow[]{st-\la}f$.

(ii) If $(f_n(x))$ is measurable $\la-a.e.$ (Definition \ref{Def1.4}) then,
\[
(f_n)\xrightarrow[]{st-\la}f\;\Rightarrow\;(f_n)\xrightarrow[]{st-a.e.}f.
\]
\end{thm}
\begin{note}\label{note1.7}
If a sequence of measurable functions $(f_n)$ converges $st-a.e.$ to some function $f$, then $f$ is measurable. Indeed, from Proposition \ref{prop1.2} (iii), it follows for $n\in\N$ that
\[
\frac{1}{n}\sum^n_{k=1}g_k(x)\overset{a.e.}{\longrightarrow}g_M(x), \ \ n\ra\infty
\]
where
\[
g_M(x)=\left\{\begin{array}{ccc}
                f(x), & \text{if} & |f(x)|\le M\\
                M, & \text{if} & f(x)> M\\
                -M, & \text{if} & f(x)<-M
              \end{array}\right., \quad g_k(x)=\left\{\begin{array}{ccc}
                f_k(x), & \text{if} & |f_k(x)|\le M\\
                M, & \text{if} & f_k(x)> M\\
                -M, & \text{if} & f_k(x)<-M.
              \end{array}\right.
\]
This implies that $g_M$ is measurable for each $M\in\N$, hence $f$ is measurable, since $g_M\overset{a.e.}{\longrightarrow}f$.
\end{note}
\begin{rem}\label{rem1.8}
If we consider the classical example of the sequence
\[
(f_n)=\big(\chi_{\big[0,\frac{1}{2}\big)},\chi_{\big[\frac{1}{2},1\big)},
\chi_{\big[0,\frac{1}{2}\big)},\ldots,\chi_{\big[\frac{i-1}{2^n},\frac{i}{2^n}\big)},
\ldots\big), \ \ i=1,2,\ld,2^n, \ \ n\in\N
\]
we easily see that $(f_n)\overset{\la}{\longrightarrow}f=0$, $(f_n)\overset{a.e.}{\not\rightarrow}f=0$ but $f_n\xrightarrow[]{st-a.e.}f=0$. Hence the following question arises:

Does it hold a weak form of the converse of Lebesgue theorem? That is, does convergence in measure, which is stronger than statistical convergence in measure imply $st-a.e.$ convergence?

In Section 2, we see that, this is not true. Also in the same section we study the notion of measurability (Definition \ref{Def1.4}) and we present an improvement of Steihaus's theorem.

Finally in Section 3 we study counterparts of the classical theorems of measure theory for statistical convergences and we find that, except Egorov's theorem, the other theorems have true corresponding versions for statistical convergences. For example the corresponding Lebesgue dominated theorem for statistical convergences is true.
\end{rem}
\section{Remarks on Steinhaus theorem}\label{sec2}
\noindent

First we see that the notion of measurability has meaning. That is, there are sequences $(a_n)$ such that the density of the sets $\{n:a_n<a\}$ exists, except for a countable number of $a\in\R$.
\begin{exm}\label{exm2.1}
Let $A_n=\{2^{n-1}\cdot k:k=1,3,5,\ld\}$, $n=1,2,3,\ld$. Then it holds that,
\[
A_n\cap A_m=\emptyset, \ \ \text{for} \ \ n\neq m \; \& \; \bigcup^\infty_{n=1}A_n=\N\;\&\; d(A_n)=\frac{1}{2^n}, \ \ n=1,2,\ld\,.
\]
(The last equality above holds since the members of each $A_n$ form an arithmetic progression with difference of successive terms equal to $2^n$).

We divide each $A_n$ into two disjoint subsets $A_{n,1}$ and $A_{n,2}$, which do not have density (e.g., $A_{n,1}$, $A_{n,2}$ are the union respectively of successive ``blocks'' in $A_n$, say $B_{1,\el}$, $B_{2,\el}$, $\el=1,2,\ld$ such that $\max B_{1,\el}<\min B_{2,\el}<\max B_{2,\el}<\min B_{1,\el+1}$ and $\dis\lim_{\el\ra\infty}\dfrac{|B_{j,\el}|}{\max B_{j,\el}}=\dfrac{1}{2^n}$ for $j=1,2$). If for each $n\in\N$, $(b_{n,\el})_{\el\in A_{n,2}}$ is an increasing sequence with first term larger than $\dfrac{1}{n+1}$ and $\dis\lim_{\el\ra\infty}b_{n,\el}=\dfrac{1}{n}$, $n=1,2,\ld$ and if for $\el\in\N$ we set,
\[
a_\el=\left\{\begin{array}{ccc}
               \frac{1}{n}, & \text{if} & \el\in A_{n,1} \\
               b_{n,\el}, & \text{if} & \el\in A_{n,2}
             \end{array}\right.,
\]
we get the following,
\begin{eqnarray}
\bigg\{\el:a_\el<\frac{1}{n}\bigg\}=\N\sm(A_1\cup\cdots\cup A_{n-1}\cup A_{1,n}), \ \ n=2,3,\ld\,.  \label{eq1}
\end{eqnarray}
\begin{eqnarray}
\{\el:a_\el<a\}=\N\sm(A_1\cup\cdots\cup A_{n-1}\cup A'_n), \ \ n=2,3,\ld,\frac{1}{n+1}<a<\frac{1}{n}.  \label{eq2}
\end{eqnarray}
The set $A'_n$ in (\ref{eq2}) differs from $A_n$ by a finite set. Hence $d(A'_n)=\frac{1}{2^n}$. Since the unions in (\ref{eq1}), (\ref{eq2}) are disjoint and the sets $A_1,A_2,\ld,A_{n-1}$, $A_n$ have densities, while the set $A_{1,n}$ does not have density, it follows that the density of the set $\{\el:a_\el<a\}$ does not exist exactly for $a=\frac{1}{2},\frac{1}{3},\ld,\frac{1}{n},\ld\,.$
\end{exm}
\begin{rems}\label{rems2.2}
(i) If a sequence of real numbers $(a_n)$ converges statistically to $\el\in\R$ then $(a_n)$ is measurable. Indeed, for each $a\in\R-\{\el\}$ it holds that
\[
d(\{n:a_n<a\})=\left\{\begin{array}{ccc}
                        1, & \text{if} & a>\el \\
                        0, & \text{if} & a<\el.
                      \end{array}\right.
\]
We note that for $a=\el$ the density of the above set may fail to exist. For example if $A$ is a subset of $\N$ which does not has density and $a_n=\el-\frac{1}{n}$ for $n\notin A$, $a_n=\el$ for $n\in A$, then $(a_n)\overset{st}{\longrightarrow}\el$ and the density of the set $\{n:a_n<\el\}$ does not exist.

(ii) If a sequence $(a_n)$ in $\R$ is measurable, then the same is true for the sequence $(|a_n|)$. Indeed, let $A=\{a\in\R:d(\{n\in\N:a_n<a\})$ does not exist$\}$ and $A'=A\cup\{a_n:n\in\N\}$. Then $A'$ is countable and for $a\notin A'\cup(-A')$, it holds that,
\begin{align*}
\{n\in\N:|a_n|<a\}&=\{n\in\N:a_n<a\}\sm\{n\in\N:a_n\le-a\} \\
&=\{n\in\N:a_n<a\}-\{n\in\N:a_n<-a\}.
\end{align*}
Hence the density of the set $\{n\in\N:|a_n|<a\}$ exist.
\end{rems}

The converse of the above implication is not true. Indeed, let $A,B$ be two disjoint subjets of $\N$, which do not have density and $A\cup B=\N$. If $(a_n)_{n\in A}$, $(b_n)_{n\in B}$ are two increasing sequences of positive real numbers, which converge to 1, we set
\[
c_n=\left\{\begin{array}{ccc}
             a_n, & \text{if} & n\in A \\
             -b_n, & \text{if} & n\in B.
           \end{array}\right.
\]
It is easy to see that the densities of the sets $\{n\in\N:c_n<a\}$ do not exist for $a\in(0,1)$, hence $(c_n)$ is not measurable, but $(|c_n|)$ is measurable since,
\[
d(\{n\in\N:|c_n|<a\})=\left\{\begin{array}{ccc}
                               1, & \text{if} & a>0 \\
                               0, & \text{if} & a\le0.
                             \end{array}\right.
\]

(iii) Let $(f_n)\xrightarrow[]{st-a.e.}f$. We may assume that $f=0$. Then,
\setcounter{equation}{0}
\begin{align}
(f_n)\xrightarrow[]{st-a.e.}f=0\Leftrightarrow&\;\fa\;\e>0:d(\{n\in\N:\,|f_n(x)|\ge\e\})=0, \ \ \la-a.e. \label{eq1} \\
\Leftrightarrow&\;\exists\;(\e_j)\in C^+_0\;\exists\;D\subseteq[0,1),\;\la(D)=0:  \label{eq2} \\
&d(\{n\in\N:|f_n(x)|\ge\e_j\})=0\;\text{for}\; x\notin D,\;j\in\N \nonumber\\
\Leftrightarrow&\;\frac{1}{n}\sum^n_{k=1}\chi_{[|f_k|\ge\e_j]}\overset{a.e.}
{\longrightarrow}0, \; n\ra\infty,\;\text{for}\; j\in\N.  \label{eq3}
\end{align}
Hence, the existence of the densities of the sets $\{n\in\N:\,|f_n(x)|<\e\}$ for all $\e>0$, $\la-ae$, is a necessary condition in order to have $st-ae$ convergence. On the other hand, in view of Remark \ref{rems2.2} (ii) above, we get a stronger form of Steinhaus' theorem, if we assume the existence of the densities of the sets $\{n\in\N:|f_n(x)|<\e_j\}$ $\la-a.e.$, $(j=1,2,\ld)$, for some $(\e_j)\in C^+_0$, or the measurability of $(|f_n(x)|)$ $\la-a.e.$. More precisely we have the following theorem.
\begin{thm}\label{thm2.3}
Suppose $(f_n)\xrightarrow[]{st-\la}f$ and that
\[
\exists\;(\e_j)\in C^+_0:d(\{n\in\N:|f_n(x)-f(x)|<\e_j\}) \;\text{exist}\; \la-a.e., \;\text{for}\; j=1,2,\ld\,.
\]
Then,
\[
\fa\;\e>0: \ \ d(\{n\in\N:|f_n(x)-f(x)|\ge\e\})=0 \ \ \la-a.e.,
\]
that is, $(f_n)\xrightarrow[]{st-a.e.}f$.
\end{thm}

The proof is similar with the proof given in
\cite{4}, with some modifications. For the reader's convenience we sketch the proof. Suppose for simplicity that $f=0$ and let $(\e_j)\in C^+_0$ such that $d(\{n\in\N:|f_n(x)|<\e_j\})$ exists $\la-a.e.$ for $j=1,2,\ld\,.$

Then
\begin{align*}
(f_n)\xrightarrow[]{st-\la}0&\Leftrightarrow\la([|f_n|\ge\e_j])\overset{st}
{\longrightarrow}0, \ \ n\ra\infty, \ \ \text{for} \ \ j=1,2,\ldots \\
&\Leftrightarrow\frac{1}{n}\sum^n_{k=1}\la([|f_k|\ge\e_j])\ra0, \ \ n\ra\infty, \ \ \text{for} \ \ j=1,2,\ld\,.
\end{align*}
(It follows by Proposition \ref{prop1.2}, (iii)).

Hence we get,
\[
\frac{1}{n}\sum^n_{k=1}\int\chi_{[|f_k|\ge\e_j]}d\la\ra0, \ \ j=1,2,\ld\,.
\]
By Fatou's Lemma, we have
\[
\lim\inf\frac{1}{n}\sum^n_{k=1}\chi_{[|f_k|\ge\e_j]}=0, \ \ \la-a.e., \ \ j=1,2,\,.
\]
Since,
\[
d(\{n\in\N:|f_n(x)|\ge\e_j\})=\lim_{n\ra0}\frac{1}{n}\sum^n_{k=1}
\chi_{[|f_k|\ge\e_j]}(x)
\]
we take that,
\[
\frac{1}{n}\sum\chi_{[|f_k|\ge\e_j]}\overset{a.e.}{\longrightarrow}0, \ \
n\ra\infty, \ \ (j=1,2,\ld).
\]
By (\ref{eq3}) above, it follows that $(f_n)\overset{a.e.}{\longrightarrow}0$.

In the next example we see that, convergence in measure, which is stronger than statistical in measure convergence, does not imply in general $st-a.e.$ convergence.
\begin{exm}\label{exm2.4}
Let $B_1=\{1,2\}$,
\[
B_n=\{n\in\N:2^1+\cdots+2^{n-1}+1\le n\le2^1+\cdots+2^{n-1}+2^n\}, \ \ n=2,3,\ld\,.
\]
We attach to each block $B_n$ $(n=1,2,\ld)$ positive integers,
\[
m_{n-1}<k^n_1<\ldots<k^{(n)}_{2^n}<m_n \ \ (m_0=0),
\]
such that
\[
\frac{k^{(n)}_1-m_{n-1}}{k^{(n)}_1}>\frac{1}{2}, \ \ \frac{k^{(n)}_{j+1}-
k^{(n)}_j}{k^{(n)}_{j+1}}>\frac{1}{2} \ \ (j=1,2,\ld,2^n-1), \ \ \frac{k^{(n)}_{2^n}}{m^n}<\frac{1}{3}.
\]
Let $I^{(n)}_j=\Big[\frac{j-1}{2^n},\frac{j}{2^n}\Big]$, $j=1,2,\ld,2^n$, $n\in\N$. We set
\[
f_k=\left\{\begin{array}{lcl}
    \chi_{I^{(n)}_1}, & \text{if} & m_{n-1}<k\le k^{(n)}_1 \\
    \chi_{I^{(n)}_j}, & \text{if} & k^{(n)}_j<k\le k^{(n)}_{j+1} \\
    0, & \text{if} & k^{(n)}_{2^n}<k<m_n
  \end{array}\right. \ \ (j=1,2,\ld,2^n-1), \ \ (n\in\N).
\]
Clearly, $(f_k)_k$ converges in measure to $f=0$.

Now, let $x\in\Ga$. Then there exist positive integers $j_1,j_2,\ld,j_n,\ld$, where $j_n\le2^n$ $(n\in\N)$, such that
\[
\{x\}=\bigcap^\infty_{n=1}I^{(n)}_{j_n}.
\]
For the corresponding increasing sequence of positive integers $k^{(1)}_{j_1}<k^{(2)}_{j_2}<\ldots<k^{(n)}_{j_n}<\ldots$ it holds that
\[
\frac{\big|\{k\le k^{(n)}_{j_n}:f_k(x)\ge\e\}\big|}{k^{(n)}_{j_n}}>
\frac{k^{(n)}_{j_n}-k^{(n)}_{j_n-1}}{k^{(n)}_{j_n}}>\frac{1}{2},
\]
where $k^{(n)}_0=m_{n-1}$, if $j_n=1$ and $\e=(0,1)$. Hence,
\begin{eqnarray}
\lim\sup\frac{|\{k\le n:f_k(x)\ge\e\}|}{n}\ge\frac{1}{2}.  \label{eq4}
\end{eqnarray}
But (\ref{eq1}), it follows that $(f_k(x))_k\overset{st-a.e.}{\not\rightarrow}0$. (In fact $f_n(x)\overset{st}{\not\rightarrow}f=0$ for all $x\in\Ga$). Also, since
\[
\frac{|\{k\le m_n:f_k(x)\ge\e\}|}{m_n}<\frac{k^{(n)}_{2^n}}{m_n}<\frac{1}{3},
\]
we have that
\begin{eqnarray}
\lim\inf\frac{|\{k\le n:f_k(x)\ge\e\}|}{n}\le\frac{1}{3}.  \label{eq5}
\end{eqnarray}
By (\ref{eq4}) and (\ref{eq5}) we see also that $(f_k(x))_k$ is not measurable for all $x\in\Ga$.
\section{The classical theorems for statistical convergences}\label{sec3}
\noindent

First, we study almost uniform convergence. We recall that a sequence is almost uniformly Cauchy, ($(f_n)$ is $a.u.-C)$, if $\fa\;\e>0$ $\exists\;D\in\Sig$, $\la(D)<\e$ such that $\big(f_n\big|_{[0,1)\sm D}\big)$ is uniformly Cauchy or equivalently,
\[
\fa\;\e'>0 \ \ \exists\;n_0\in\N:\sup_{x\notin D}|f_n(x)-f_{n_0}(x)|<\e' \ \ \text{for} \ \ n\ge n_0.
\]
Moreover, it is well known that,
\setcounter{equation}{0}
\begin{eqnarray}
(f_n)_n \ \ \text{is} \ \ a.u-C\Leftrightarrow\exists\;f\in L^0:f_n\xrightarrow[]{a.u}f.  \label{eq1}
\end{eqnarray}
These notions are generalized naturally for statistical convergences.
\end{exm}
\begin{Def}\label{Def3.1}
(a) We say that, $(f_n)$ is statistically almost uniformly Cauchy $(st-a.u-C)$, if
\begin{align*}
\fa\;\e>0\ \ \exists\; D\in\Sig,\ \ &\la(D)<\e\ \ \fa\;\e'>0\ \ \exists\;n_0\in\N:\\
&d\Big(\Big\{n\in\N:\sup_{x\notin D}|f_n(x)-f_{n_0}(x)|\ge\e'\Big\}\Big)=0.
\end{align*}

(b) We say that $(f_n)$ converges statistically almost uniformly to $f$, if
\[
\fa\;\e>0\ \ \exists\; D\in\Sig,\ \ \la(D)<\e\ \ \fa\;\e'>0\ \
d\Big(\Big\{n\in\N:\sup_{x\notin D}|f_n(x)-f(x)|\ge\e'\Big\}\Big)=0.
\]
\end{Def}

For the proof of the next theorem we will need the following lemma, which asserts that the $n_0\in\N$ in Definition \ref{Def3.1} (a) can be chosen arbitrarily large in any set of density 1.
\begin{lem}\label{lem3.2}
The following are equivalent

(i) $(f_n)$ is $st-a.u-C$.

(ii) $\fa\;B\subseteq\N$, $d(B)=1$ $\fa\;\e>0$ $\exists\;D\in\Sig$, $\la(D)<\e$, $\fa\;\e'>0$ $\fa\;N\in B$ $\exists\;n_1\in B$, $n_1\ge N$: $d\Big(\Big\{n\in B:\dis\sup_{x\notin D}|f_n(x)-f_{n_1}(x)|<\e'\Big\}\Big)=1$.
\end{lem}
\begin{Proof}
Suppose that $(f_n)$ is $st-a.u.-C$ and let $B\subseteq\N$ with $d(B)=1$, $\e>0$, $D\in\Sig$ with $\la(D)<\e$, $\e'>0$ and $N\in B$. Then, by hypothesis, there exists $n_0\in\N$ such that $d(B_0)=1$, where $B_0=\Big\{n\in\N:\dis\sup_{x\notin D}|f_n(x)-f_{n_0}(x)|<\frac{\e'}{2}\Big\}$.

If $n_1\in B_0\cap B$ with $n_1\ge n_0,N$ and $n\in B\cap B_0$ with $n\ge n_1$, then
\[
\sup_{x\notin D}|f_n(x)-f_{n_1}(x)|\le\sup_{x\notin D}|f_n(x)-f_{n_0}(x)|+\sup_{x\notin D}|f_{n_1}(x)-f_{n_0}(x)|<\e'.
\]
Hence,
\[
\{n\in B\cap B_0:n\ge n_1\}\subseteq\Big\{n\in B:\sup_{x\notin D}|f_n(x)-f_{n_1}(x)|<\e'\Big\}.
\]
Since the density of the first set in the above inclusion is 1, the results follows.

The converse implication obviously holds.
\end{Proof}
\begin{thm}\label{thm3.3}
The following are equivalent
\begin{enumerate}
\item[(I)] $(f_n)_n$ is $st-a.u.-C$
\item[(II)] $\exists\;B\subseteq\N$, $d(N)=1$: $(f_n)_{n\in B}$ is $a.u.-C$
\item[(III)] $\exists\;B\subseteq N$, $d(B)=1$ $\exists\;f\in L^0$: $(f_n)_{n\in B}\overset{a.u}{\longrightarrow}f$
\item[(IV)] $\exists\;f\in L^0$: $(f_n)_n\xrightarrow[]{st-a.u}f$.
\end{enumerate}
\end{thm}
\begin{Proof}
(I)\;$\Rightarrow$\;(II)

Suppose that $(f_n)_n$ is $st-a.u.-C$. We construct by induction an increasing sequence $(n_k)_k$ in $\N$, a sequence $(C_k)_k$ in $\Sig$ and a decreasing sequence $(B_k)_k$ of subsets of $\N$ such that

1. $\mi(C_k)<\frac{1}{2^k}$, $k=1,2,\ld$

2. $d(B_k)=1$ and $B_k=\Big\{n\in B_{k-1}:\dis\sup_{x\notin C_k}|f_n(x)-f_{n_k}(x)|<\frac{1}{k}\Big\}$, for $k=1,2,\ld$ $(B_0=\N)$.

3. $n_k\in B_k$, $k=1,2,\ld$\smallskip

4. $\dfrac{|\{n\in B_k:m\le n\}|}{n}>1-\dfrac{1}{k+1}$, for $n\ge n_{k+1}$, $k=1,2,\ld\,.$\medskip\\
\noindent
{\bf Step 1.} Let $B_0=\N$. By hypothesis there exist $n_1\in B_0$, $C_1\in\Sig$ with $\mi(C_1)<\frac{1}{2}$ such that $d(B_1)=1$, where $B_1=\Big\{n\in B_0:\dis\sup_{x\notin C_1}|f_n(x)-f_{n_1}(x)|<\frac{1}{1}\Big\}$.

Since $d(B_1)=1$, we can find $n'_2\in B_1$, $n'_2>n_1$ such that $|\{m\in B_1:m\le n\}|\big/n>1-\frac{1}{1+1}$ for $n\ge n'_2$.

By Lemma \ref{lem3.2} it follows that
\[
\exists\;C_2\in\Sig, \ \ \mi(C_2)<\frac{1}{2^2}\ \ \exists\;n_2\in B_1, \ \ n_2>n'_2:d(B_2)=1, \ \ \text{where}
\]
\[
B_2=\Big\{n\in B_1:\sup_{x\notin C_2}|f_n(x)-f_{n_2}(x)|<\frac{1}{2}\Big\}.
\]
Hence 1,2,3,4 are satisfied for $k=1$ and simultaneously 1, 2, 3 for $k=2$.\medskip \\
\noindent
{\bf Step $\bbb{k}$.} Suppose we have defined $n_1<n_2<\ldots<n_k$, $C_1,C_2,\ld,C_k\in\Sig$ and $B_0=\N\supseteq B_1\supseteq B_2\supseteq\cdots\supseteq B_k$ such that 1,2,3,4 hold for $1,2,\ld,k-1$ and 1, 2, 3 hold for $1,2,\ld,k$. Since $d(B_k)=1$, we can find $n'_{k+1}\in B_k$, $n'_{k+1}>n_k$, such that $|\{m\in B_k:m\le n\}|\big/n>1-\frac{1}{k+1}$ for $n\ge n'_{k+1}$. Again by Lemma \ref{lem3.2} it follows that
\[
\exists\;C_{k+1}\in\Sig, \ \ \mi(C_{k+1})<\frac{1}{2^{k+1}} \ \ \exists\;n_{k+1}\in B_k, \ \ n_{k+1}>n'_{k+1}:d(B_{k+1})=1, \ \ \text{where}
\]
\[
B_{k+1}=\Big\{n\in B_k:\sup_{x\notin C_{k+1}}|f_n(x)-f_{n_{k+1}}(x)|<\frac{1}{k+1}\Big\}.
\]
So the induction processes is completed.

We set,
\[
B=\{n\in B_1:n\le n_3\}\cup\{n\in B_2:n\leqq n_4\}\cup\ld,
\]
then for $n_{k+1}\le n<n_{k+2}$ it holds that
\[
\{m\in B:m\le n\}\supseteq\{m\in B_k:m\le n\}.
\]
Hence,
\[
\frac{|\{m\in B:m\le n\}|}{n}\ge1-\frac{1}{k+1}, \ \ \text{for} \ \ n_{k+1}\le n,
\]
which implies that $d(B)=1$.

Now, let $\e>0$ and $k\in\N$ such that $\ssum^\infty_{\el=k}\frac{1}{2^\el}<\e$. If $D=\bigcup\limits^\infty_{\el=k}C_\el$, then for any $\e'>0$ and $k'\in\N$, $\frac{1}{k'}<\e'$, $k'\ge k$, it follows that
\[
\sup_{x\notin D}|f_n(x)-f_{n_{k'}}|\le\sup_{x\notin C_{k'}}|f_n(x)-f_{n_{k'}}|<\frac{1}{k'}<\e,
\]
for all $n\in B$, $n\ge n_{k'}$ (As $\{n\in B:n\ge n_{k'}\}\subseteq B_{k'}$ and the last inequality above holds, by 2, for all $n\ge n_{k'}$). This means that the sequence $(f_n)_{n\in B}$ is $a.u.-C$ and the proof is complete.
\end{Proof}
\noindent
(II)\;$\Rightarrow$\;(III)

It follows by (\ref{eq1}) \medskip \\
\noindent
(III)\;$\Rightarrow$\;(IV)

It follows at once that Definition \ref{Def3.1} (b) is satisfied.\medskip \\
\noindent
(IV)\;$\Rightarrow$\;(I)

Let $\e>0$. Then by hypothesis there exists $D\in\Sig$, $\mi(D)<e$ such that for any $\e'>0$ it holds that $d(M)=1$, where $M=\Big\{n\in\N:\dis\sup_{x\notin D}|f_n(x)-f(x)|<\frac{\e'}{2}\Big\}$. We fix $n_0\in M$. Then for all $n\in M$ it holds that
\[
\sup_{x\notin D}|f_n(x)-f_{n_0}(x)|\le\sup_{x\notin D}|f_n(x)-f(x)|+
\sup_{x\notin D}|f_{n_0}(x)-f(x)|<\e'.
\]
Hence,
\[
M\subseteq\Big\{n\in\N:\sup_{n\notin D}|f_n(x)-f_{n_0}(x)|<\e'\Big\}
\]
which implies that $(f_n)_n$ is $st-a.u.-C$.
\begin{rem}\label{rem3.4}
We note that Theorem \ref{thm3.3} holds for arbitrary measure spaces, not necessarily finite, as we consider in this paper.
\end{rem}

The classical Riesz theorem (see \cite{6}) asserts that, a sequence $(f_n)$ converges in measure to some $f\in L^0$, if and only if, $(f_n)$ is Cauchy in measure (that is, $\fa\;\e>0$ $\fa\;\de>0$, $\exists\;n_0\in\N$, such that $\la([|f_n-f_{n_0}|\ge\e])<\de$ for all $n\ge n_0$).

On the other hand, the following facts are well known:
\begin{enumerate}
\item[$\bullet$] A sequence $(f_n)$ converges in measure of $f$, if and only if, $(f_n)$ converges to $f$ with respect to the following metric $\rho$ on the space $L^0$:
    \[
    \rho(f,g)=\inf\{\e+\la([|g-g|\ge\e]):\e>0\}
    \]
(see \cite{3}).
\item[$\bullet$] A sequence $(f_n)$ is Cauchy in measure, if and only if, $(f_n)$ is Cauchy sequence with respect to the metric $\rho$.
\end{enumerate}

It is not hard to see that,
\begin{eqnarray}
(f_n)\xrightarrow[]{st-\la}f\Leftrightarrow(f_n)\xrightarrow[]{st-\rho}f \label{eq2}
\end{eqnarray}
(that is $(f_n)$ converges statistically to $f$ in the metric space $(L^0,\rho)$. See also Remark \ref{rem1.3}). Indeed, if $(f_n)\xrightarrow[]{st-\la}f$ then by Definition \ref{Def1.5} (b), it follows that, there are disjoint finite subsets $B_k$ of $\N$, $k=1,2,\ld$, with max $B_k<\min B_{k+1}$ for all $k\in\N$ such that
\[
\la\bigg(\bigg[|f_n-f|\ge\frac{1}{k}\bigg]\bigg)<\frac{1}{k} \ \ \text{for} \ \ n\in B_k \ \ \& \ \ \frac{|B_k|}{\max B_k}>1-\frac{1}{k}.
\]
If we set $B=\bigcup\limits^\infty_{k=1}B_k$ then, we easily see that,
\[
d(B)=1 \ \ \& \ \ (f_n)_{n\in B}\overset{\rho}{\longrightarrow}f.
\]
Hence
\[
(f_n)\xrightarrow[]{st-\rho}f.
\]
Conversely if $(f_n)\xrightarrow[]{st-\rho}f$, then by definition of the metric $\rho$, it follows at once that, $(f_n)\xrightarrow[]{st-\la}f$.

Also, if we define $(f_n)$ to be statistically Cauchy in measure (in symbols $(f_n)$ is $st-C-\la$), if and only if, $\fa\;\e>0$ $\fa\;\de>0$ $\exists\;n_0\in\N$:
\[
d\big(\big\{n\in\N:\la\big([|f_n-f_{n_0}|\ge\e]\big)\ge\de\big\}\big)=0,
\]
then similarly as (\ref{eq2}) above we get,
\[
(f_n) \ \ \text{is} \ \ st-C-\la\Leftrightarrow(f_n) \ \ \text{is} \ \ st-C \ \ \text{in} \ \ (L^0,\rho).
\]
(See Definition \ref{Defs1.1} (b) and Remark \ref{rem1.3}).

Hence, the proof of the corresponding version of Riesz theorem for statistical converges follows from Proposition \ref{prop1.2} (ii) and Remark \ref{rem1.3}:
\begin{thm}\label{thm3.4}
The following are equivalent

(i) $(f_n)$ is $st-C-\la$

(ii) $\exists\;f\in L^0:(f_n)\xrightarrow[]{st-\la}f$.
\end{thm}

Now, we turn to Egorov's theorem. It is known that, this theorem is not true for statistical convergences, that is, $st-a.e.$ convergence does not imply in general $st-a.u$ convergence. This fact is contained in \cite{1} (see [1] \S\,3, Theorem 11 and Example 13). Here we present a much simpler example than that of \cite{1}, which assures that Egorov's theorem is not true for statistical converges.
\begin{exm}\label{exm3.4}
We consider the sequence $(f_k)$, where $f_k=\chi_{A_k}=\chi_{\big[\frac{j-1}{2^n},\frac{j}{2^n}\big]}$, if $k=2^1+2^2+\cdots+2^{n-1}+j$, $j\in\{1,2,\ld,2^n\}$. As we saw in Remark \ref{rem1.8} it holds that
\[
(f_k)\xrightarrow[]{st-a.e.}f=0.
\]
{\em Assertion}. If $K=\{k_1<k_2<\ldots<k_m<\ldots\}\subseteq\N$, with $d(K)=1$, then $x\in A_{k_m}$ for infinitely many $m\in\N$ \, $\la-a.e.$\vspace*{0.2cm} \\
\noindent
{\em Proof of Assertion}. We set,
\[
K'=\N\backslash K \ \ \text{and} \ \ C_n=\bigcup_{k\in K'\cap B_n}A_k
\]
where $B_n=\{2^1+2^2+\cdots+2^{n-1}+1,\ld,2^1+2^2+\cdots+2^{n-1}+2^n\}$, $n\in\N$.

It suffices to show that
\[
\la(\lim\inf C_n)=\int\lim\inf\chi_{C_n}d\la=0.
\]
By Fatou's Lemma we have:
\begin{eqnarray}
\int\lim\inf\chi_{C_n}d\la\le\lim\inf\int\chi_{C_n}d\la.  \label{eq3}
\end{eqnarray}
But,
\begin{eqnarray}
\int\chi_{C_n}d\la=\frac{|k\in\N:k\in K'\cap B_n|}{2^n}  \label{eq4}
\end{eqnarray}
and
\begin{eqnarray}
\frac{|\{k\in\N:k\in K'\cap B_n\}|}{2^1+2^2+\cdots+2^n}=
\frac{|\{k\in\N:k\in K'\cap B_n\}|}{2^n}\cdot\frac{2^n}{2^1+\cdots+2^n}.
\label{eq5}
\end{eqnarray}
Since
\[
\frac{|\{k\in\N:k\in K'\cap B_n\}|}{2^1+2^2+\cdots+2^n}\le
\frac{|\{k\in K':k\le 2^n\}|}{2^1+2^2+\cdots+2^n}
\]
and
\[
\frac{2^n}{2+2^n+\cdots+2^n}=\frac{1}{\frac{1}{2^{n-1}}+\frac{1}{2^{n-2}}+\cdots+1}
\ra\frac{1}{2}, \ \ \text{as} \ \ n\ra\infty
\]
and $d(K')=0$, it follows from (\ref{eq4}) and (\ref{eq5}) that $\dis\lim_{n\ra0}\int\chi_{C_n}d\la=0$, which in view of (\ref{eq3}) completes the proof of the assertion.

Now from the assertion it follows that
\[
(f_{k_n})\overset{a.e.}{\not\rightarrow}f=0 \ \ \text{if} \ \ d(\{k_1<\ldots<k_n<\ldots\})=1
\]
which implies, by Theorem \ref{thm3.3}, that $(f_n)\not\xrightarrow[]{st-a.u.}f=0$.
\end{exm}
\begin{prop}\label{prop3.7} Suppose that $(f_n)\xrightarrow[]{st-a.e.}f$ and that there exists $g\in L^1(\la)$ with $|f_n(x)|\le g(x)\mi-a.e.$ for $n=1,2,\ld\,.$ Then:
\[
\int f_nd\la\overset{st}{\longrightarrow}\int fd\la.
\]
\end{prop}
\begin{Proof}
It is easy to see that if $(f_n)\xrightarrow[]{st-a.e.}f$, then
\[
(f^+_n)\xrightarrow[]{st-a.e.}f^+ \ \ \text{and} \ \ (f^-_n)\xrightarrow[]{st-a.e.}f^-,
\]
where $f^+_n$, $f^-_n$ are the positive and negative parts of $f_n$, $n=1,2,\ld\,.
$

Also, if
\[
\int f^+_nd\la\overset{st}{\longrightarrow}\int f^+d\la \ \ \text{and} \ \ \int f^-_nd\la\overset{st}{\longrightarrow}\int f^-_nd\la
\]
then $\int fd\la\overset{st}{\longrightarrow}\int fd\la$, hence it is enough to assume that $f_n(x)\ge0$ $\la-a.e.$ for all $n\in\N$ and $f=0$.

Now, by hypothesis and Proposition \ref{prop1.2} (iii) we get that
\[
\frac{1}{n}\sum^n_{k=1}f_k(x)\ra0 \ \ a.e. \ \ \text{for} \ \ x\in[0,1).
\]
Hence, by Lebesgue's dominated theorem we have that
\[
\frac{1}{n}\sum^n_{k=1}\int f_kd\la\ra0,
\]
which by Proposition \ref{prop1.2} (iv) implies
\[
\int f_nd\la\overset{st}{\longrightarrow}\int fd\la.
\]
\end{Proof}
\begin{exm}\label{exm3.8}
Let $A\in\Sig$ with $\la(A)>0$ and $K=\{k_1<\ldots<k_n<\ld\}\subseteq\N$ with $d(K)=1$. If we set
\[
f_n=\left\{\begin{array}{lll}
             \chi_A, & \text{for}& n\in K \\
             0, & \text{otherwise}. &
           \end{array}\right.
\]
Then
\[
(f_n)_n\xrightarrow[]{st-a.e.}\chi_A=f \ \ \text{and} \ \ \int f_nd\la\overset{st}{\longrightarrow}\int fd\la.
\]
But
\[
\int f_nd\la\not\ra\int fd\la.
\]
Hence, the above result is the best possible regarding convergence of integrals in case of statistical convergences.
\end{exm}
\begin{rem}\label{rem3.9}
Apparently Proposition \ref{prop3.7} holds for arbitrary measure spaces, not necessarily finite.
\end{rem}\medskip
\noindent
{\bf Acknowledgements}

I thank Professor M. Kolountzakis and N. Papanastassiou for their suggestions and remarks.
\vspace*{0.6cm}
C. Papachristodoulos \vspace*{0.1cm} \\
Department of Mathematics  \\
University of Crete  \\
KNOSSOS AV. 71409  \\
HERAKLION - CRETE - GREECE  \\
e-mail: papach@math.uoc.gr


\begin{thebibliography}{99}
%
\bibitem{1} Balcerzak K., Dens K., Kemisarski A., Statistical convergence and ideal convergence of sequences of functions, J. of Math. Analysis and Its Appl., 328(2007), 715-729.
%
\bibitem{2} Connor J. S., The Statistical and strong $p$-Cesaro convergence of sequences, Analysis 8, 1988, p. 47-63.
%
\bibitem{3} Dunford N. and Schwartz J., Linear Operators Part I, Intersience Publishers, Inc., N.Y.
%
\bibitem{4} Fast H., Sur la convergence statistique, Coll. Math. 2, 1951, p. 241-244.
%
\bibitem{5} Fridy J. A., On statistical convergence, Analysis 5 (1985), 301-313.
%
\bibitem{6} Halmos P., Measure Theory, Springer-Verlag.
%
\bibitem{7} \u{S}al\'{a}t T., On Statistically convergent sequences of real numbers, Math Slovaca 30, 1980, No 2, p. 139-150.
%
\bibitem{8} H. Steinhaus, sur la convergence ordinarie et la convergence asymptotic, Colloq. math. 2 (1951) 73-74.
%
\bibitem{7} A. Zygmund, trig. series, 3nd Edition. Cambridge Math. Library.
\end{thebibliography}
\end{document}